\documentclass{dcds}
\usepackage{relsize}

\usepackage{amsmath,amssymb}
\usepackage[utf8]{inputenc}
\usepackage[english]{babel}
\usepackage{csquotes}
  \usepackage{paralist}
 \usepackage{psfrag, subfigure}
  \usepackage{graphics}
  \usepackage{epsfig}
\usepackage{graphicx}
\usepackage{epstopdf}
 \usepackage[colorlinks=true]{hyperref}
\hypersetup{urlcolor=blue, citecolor=red}
\usepackage{tikz-cd}

\def\currentvolume{}
 \def\currentissue{}
  \def\currentyear{}
   \def\currentmonth{}
    \def\ppages{X--XX}
     \def\DOI{}

\newtheorem{thm}{Theorem}[section]

\newtheorem{lema}{Lemma}[section]

\newtheorem{prop}{Proposition}[section]

\theoremstyle{definition}

\newtheorem{rk}{Remark}[section]

\newcommand{\pf}{{\flushleft{\bf Proof: }}}
\newcommand{\N}{\mbox{$\mathbb{N}$}}
\newcommand{\T}{\mbox{$\mathbb{T}$}}

\newcommand{\R}{\mbox{$I\!\!R$}}

\numberwithin{equation}{section}

\begin{document}

\title[ Robust transitivity of singular endomorphisms ]{ High dimensional robustly singular maps that are $C^1$ robustly transitive.}

\author[Juan C. Morelli]{}


\email{jmorelli@fing.edu.uy }

\maketitle

\centerline{\scshape  Juan Carlos Morelli$^*$}
\medskip
{\footnotesize
 \centerline{Universidad de La Rep\'ublica. Facultad de Ingenieria. IMERL.}
   \centerline{ Julio Herrera y Reissig 565. C.P. 11300.}
   \centerline{ Montevideo, Uruguay.}}

\bigskip

 \centerline{(Communicated by )}

\begin{abstract}

 We improve previous results by exhibiting a construction that contains all known examples. A sufficient condition for the existence of robustly transitive maps displaying singularities on a certain large class of compact manifolds is given.

\end{abstract}

\section{Introduction}

Let $X$ be a real (riemannian) manifold and an induced discrete dynamical system over $X$ by a continuous map $f: X \to X$. Many of the features of interest displayed by the dynamical system depend on the properties displayed by the map $f$, and whichever property $f$ possesses is said to be \textit{robust} if all maps sufficiently close to $f$ also posses that same property. Our article addresses robustness of transitivity, meaning by \textit{transitive} the existence of a forward dense orbit for some point $x$ in $X$.\\
Robust transitivity has been an active topic for some decades now. The first results showed that for \textsl{diffeomorphisms} it is needed to have at least weak hyperbolicity for $C^1$ robust transitivity (see \cite{m2}, \cite{bdp}), while \textit{Anosov on $\T^2$} is an equivalent condition regarding surfaces. In the setting of \textsl{regular endomorphisms} (empty critical set), weak forms of hyperbolicity are also needed for robust transitivity. And some sufficient conditions are provided for the $n$-torus as phase space (see \cite{and}, \cite{p}, \cite{lp}).  \\
Robust transitivity of \textsl{singular endomorphisms} has been the least studied setting, only in 2013 and 2016 were exhibited the first such examples (see \cite{br}, \cite{ilp}) on the $2$-torus. The setting became active afterwards, in 2019 was shown that weak hyperbolicity is needed for $C^1$ robust transitivity of singular maps on any surface (see \cite{lr1}, \cite{lr2}). The higher dimensional setting has been approached even less, with only examples exhibited on $\T^2$ now extended to $\T^n$ (see \cite{mo}, \cite{mo2}).\\
The aim of our article is to provide with a construction that contains all of the previous examples of $C^1$ robustly transitive maps displaying critical points as particular cases. We give some definitions now in order to clearly state our main result.\\ The \textit{critical set} of the map $f$ is $S_f:=\{ x \in X : det(D_xf)=0\}$, the set of points where the determinant of its jacobian matrix is null. The map $f$ is said to be \textit{singular} if its critical set is nonempty, and it is said to be \textit{robustly singular} if there exists a neighborhood $\mathcal{U}_{f}$ of $f$ in the $C^1$ topology such that the critical set $S_g$ is nonempty for all $ g \in \mathcal{U}_{f}$. We say that the map $f$ is $C^k$ \textit{robustly transitive} if there exists a neighborhood $\mathcal{U}_{f}$ of $f$ in the $C^k$ topology such that $g$ is transitive for all $ g \in \mathcal{U}_{f}$. The main result provided by our work can be stated as follows:

\begin{thm}\label{main}
  Let $M_1$ be a compact boundaryless real manifold supporting an expanding\footnote{We mean a map where all vectors of the tangent spaces are forwardly uniformly expanded by the differential map. A precise definition given by \cite{shub} is provided in the next section.} map and $M_2$ a compact boundaryless real manifold.
  Then, there exists an endomorphism of $M_1 \times M_2$ that is $C^1$ robustly transitive and robustly singular.
\end{thm}

  Theorem \ref{main} gives a sufficient condition for the existence of such a map on a large class of compact manifolds obtained as products. Many straightforward examples are obtained from it: \textit{Every elliptical ring cyclide supports a robustly singular endomorphism that is $C^1$ robustly transitive.} Or \textit{for every $n \geq 2$, the $n$-torus supports a robustly singular endomorphisms that is $C^1$ robustly transitive}. Also \textit{any knot obtained from tori do too}, etcetera. Even when Theorem \ref{main} is quite a general result containing a large number of examples, many open questions still remain. \textit{Are there other types of manifolds supporting this kind of endomorphism?} \textit{Can we relax the expanding hypothesis to (weakly) hyperbolic?} \textit{Can this construction be set in non-compact environments?} \textit{Is this construction suitable for fiber bundles?} Just to mention some of them.
\subsection{Sketch of the Construction/Contributions of the paper.}
    We start with the construction of a map $f$ supported on $M_1 \times M_2$ such that its dynamics are determined by the local behavior around a ''small'' (so called \textit{blending}) region that gives a transitive map. The construction is carried on in a $C^1$ robust way. Since this kind of constructions allow us to determine global behavior arising from local behavior, they give us a lot of freedom to perform perturbations away from the blending regions without destroying the dynamical characteristics. In this fashion we introduce critical points artificially, away from the blending region with standard surgical procedures and in a robust way. A map satisfying the claim at the thesis of Theorem \ref{main} arises.\\
    The general idea of the construction may appear unoriginal at first since robust transitivity is obtained in this way frequently. Moreover, the examples displaying singularities provided by \cite{ilp}, \cite{ip} and \cite{mo2} are constructed in this fashion; but may the reader recall that one of our purposes is to provide with a construction that overviews all previous examples. Even so, we believe our work presents a fresh approach on the topic for two reasons: it takes on the high dimensional context and it provides with a new technique for robust transitivity of singular maps. \\ All known examples prior to our paper obtain robust transitivity relying on the existence of a field of unstable cones. We dismiss this condition and obtain the property from a Cantor set in the phase space carrying special properties, usually referred to as blending properties. \\
    To finish, we point out the main obstacle for $C^1$ robust transitivity in the singular setting. May the reader observe that mostly ever a map $f$ admits critical points, for every $C^1$ neighborhood $\mathcal{U}_f$ of $f$, there exists a map $g$ in $\mathcal{U}_f$ such that for some open set $U$ of $X$, the image $g(U)$ has empty interior (for proofs see \cite[Lemma 4.0.2]{ip} and \cite[Lemma 4.5]{mo2}). This is evidently a huge obstruction for robust transitivity since if this $g(U)$ happened to be contained in a meager invariant set (which is not unusual, i.e. a sumbanifold) the transitivity property would then be lost.

\section{Preliminaries}

Let $X$ be a differentiable manifold of dimension $m$ and a differentiable endomorphism $f:X \rightarrow X$. The \textit{orbit} of $x \in X$ is $\mathcal{O}(x):=\{ f^n(x) , n \in \N \}$, and $f$ is \textit{transitive} if there exists a point $x \in X$ such that the closure of its orbit $ \overline{{\mathcal{O}(x)}}=X$. The following proposition is of major practical use:

\begin{prop}\label{equi}
  If $f$ is continuous, the following are equivalent:
  \begin{enumerate}
    \item  $f$ is transitive.
    \item  For all $U, V$ open sets in $X$, there exists $n \in \N$ such that $ f^n(U) \cap V \neq \emptyset$.
    \item  There exists a residual set $R$ (countable intersection of open and dense sets) such that all points $  x \in R$ satisfy $\overline{\mathcal{O}(x)}=X$.
  \end{enumerate}
\end{prop}

 \subsection{Blenders.}

 A brief overview of the concept of a \textit{blender} is given now. In most situations it is easy to think of blenders as higher dimensional horseshoes, or as sets exhibiting the dynamics of a Smale's horseshoe. Blenders force the robust intersection of topologically 'thin' sets, giving rise to rich dynamics.\\ According to \cite{bcdw}, \begin{displayquote}{\textit{''A blender is a compact hyperbolic set whose unstable set has dimension strictly less than one would predict by looking at its intersection with families of submanifolds''.}}\end{displayquote}

\begin{figure}[ht]
\begin{center}
\includegraphics[scale=0.4]{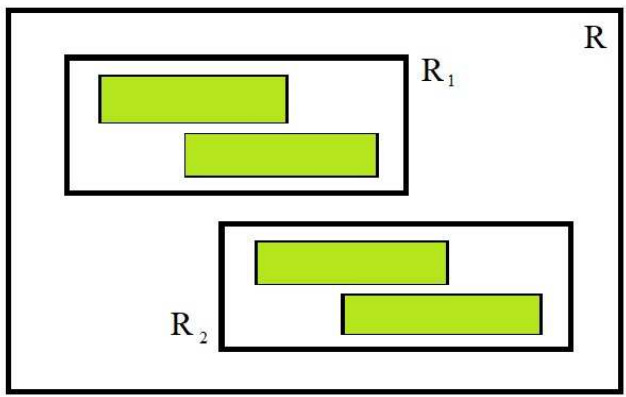}
\caption{A protoblender over $R$. Darker is $f^{-1}(R)$.}\label{protoblender}
\end{center}
\end{figure}

 They also provide with a prototipical example of a blender: Let $R$ be a rectangle with two rectangles $R_1$ and $R_2$ lying inside, horizontally, and such that their projections onto the base of $R$ overlap (Figure \ref{protoblender}). Consider now a diffeomorphism $f$ such that $f(R_1)=f(R_2)=R$. Then, $\Omega= \bigcap_{n \in \N}f^{-n}(R)$ gives rise to a blender (Cantor) set for $f$. Observe that $f$ admits a fixed point inside each of $R_1$ and $R_2$, and that all vertical segments between the projection of these points intersect $\Omega$ (this is due to the overlapping of the projections of $R_1$ and $R_2$ which holds at every preiteration). Observe as well that this construction is robust in two senses: on the one hand, $f$ can be slightly perturbed with persistance of the property. And on the other, the vertical segment can also be slightly perturbed and still intersect $\Omega$.

Notice that $\Omega$ is a fractal object with topological dimension zero. Nonetheless, every close-to-vertical line in between the fixed points of $f$ inside $R_1 \cup R_2$ intersects $\Omega$; hence, one would expect $\Omega$ to be at least of topological dimension one. This is the characteristical trait of blender sets.\\
To finish with the preliminaries regarding blenders, their importance lies in the fact that they are a magnificent tool for producing rich dynamics, particularly robustly transitive dynamics in the absence of hyperbolicity. For more insight on blenders and its applications the reader may go to \cite{bd} or \cite{bcdw}.

\subsection{Iterated Function Systems.}

Let $\mathcal{F},\mathcal{G}$ be two families of diffeomorphisms of $M$. Denote by $\mathcal{F} \circ \mathcal{G}:= \{f \circ g  / \quad f \in \mathcal{F}, g \in \mathcal{G}\}$; and for $k \in \N$ denote $\mathcal{F}^0=\{ Id_M\}$ and $\mathcal{F}^{k+1}=\mathcal{F}^{k} \circ \mathcal{F}$. Then, the set $\bigcup_{k=0}^{\infty}\mathcal{F}^k$ has a semigroup structure that is denoted by $\langle \mathcal{F}\rangle^+$ and said to be generated by $\mathcal{F}$. The action of the semigroup $\langle \mathcal{F}\rangle^+$ on $M$ is called the \textit{iterated function system} associated with $\mathcal{F}$, we denote it by IFS$(\mathcal{F})$; and for $x \in M$, the \textit{orbit} of $x$ by the action of the semigroup $\langle \mathcal{F}\rangle^+$ is $ \langle \mathcal{F}\rangle^+(x)=\{f(x), f \in \langle \mathcal{F}\rangle^+ \}$. A sequence $\{ x_n \}_{ n \in \N }$ is a branch of an orbit of IFS$(\mathcal{F})$ if for every $n \in \N$ there exists $f_n \in \langle \mathcal{F}\rangle^+$ such that $f_n(x_n)=x_{n+1}$.\\ An IFS$(\mathcal{F})$ is \textit{minimal} if for every $x \in M$ the orbit $ \langle \mathcal{F}\rangle^+(x)$ has a branch that is dense on $M$; and is \textit{$C^r$ robustly minimal} if for every family $\mathcal{\hat{F}}$ of $C^r$ perturbations of $\mathcal{F}$ and every $x \in M$ the orbit $\langle \mathcal{\hat{F}}\rangle^+(x)$ has a branch that is dense on $M$.\\
The next theorem is crucial for the construction carried on in the article. For the proof the reader may go to \cite[Theorem B]{hn}.

\begin{thm}\label{SRM}
  Every boundaryless compact manifold admits a pair of diffeomorphisms $\{g_1,g_2\}$ that generate a $C^1$ robustly minimal IFS. In particular, the diffeomorphism $g_1$ admits a unique attracting fixed point $a_1$ and a unique repelling fixed point $r_1$.
\end{thm}

Having stated all the preliminary facts needed to construct the map claimed to exist by Theorem \ref{main}, we proceed to it now in two steps. Define first an endomorphism $f$ of $M_1 \times M_2$ that is $C^1$ robustly transitive and then perturb $f$ to a map $A$ satisfying the claim at the thesis of Theorem \ref{main}.

\section{A robustly transitive map $f$ of $M_1 \times M_2$.}
\subsection{Construction of $f$.}
Starting from an expanding map $F$ supported on $M_1$ we construct a perturbed map of $F \times Id_{M_2}$ displaying a point with a neighborhood whose forward and backward iterates are dense on $M_1 \times M_2$.\\

 Let $M_1$ and $M_2$ be compact real manifolds of dimension $m_1$ and $m_2$ respectively, endowed with riemannian metrics (without loss of generality consider the metrics to be euclidean). Let ${f_0}$ be an expanding endomorphism supported on $M_1$. Recall that, in the sense of \cite{shub}, ${f_0} \in C^1(M_1)$ is \textit{expanding} if there exist constants $c>0$ and $k>1$ such that for every point $x \in M_1$ and every tangent vector $v \in T_x M_1$ holds that $||D_xf_0^n(v)||> ck^n||v||$, for all $n \in \N^+$. Recall as well that every expanding map admits a fixed point $p$ \cite[Theorem 1 (a)]{shub}.\\
 Let $U$ be a connected open neighborhood of $p$, fixed point of $f_0$ in $M_1$ and let $V$ be an open connected set in $M_1$ disjoint from $U$. Define for all subsets $X$ of $M_1$ the set $X_r := \bigcup_{x \in X}B_{(x,r)}$ and fix $\varepsilon >0$ such that $U_\varepsilon \cap V_\varepsilon = \emptyset$. \\
 Since $f_0$ is expanding, by \cite[Theorem 1 (c)]{shub} there exist $n_1,n_2 \in \N$ satisfying $f_0^{n_1}(U)=f_0^{n_2}(V)=M_1$ . Define a new map $F: M_1 \rightarrow M_1$ by $F:=f_0^{\max\{n_1,n_2\}}$ so we have constructed an expanding endomorphism $F$ supported on $M_1$ satisfying $F(U)=F(V)=M_1$ and $F(p)=p$. Observe that there is no loss in considering that the determinant \textit{det}$(D_xF)$ is positive for all $x \in M_1$ (take a power of $F$ if needed).

 \begin{rk}\label{preimagenes}
 \begin{enumerate}\
   \item The sets $U$ and $V$ define a \textit{proto-blender-like} structure for $F$.
  \item The set $C=\bigcap_{n \in \N} F^{-n}(U \cup V)$ is a \textit{blending-region-like} Cantor set.
 \end{enumerate}

 \end{rk}

 Define now a smooth ($C^\infty$) bump function $u: M_1 \rightarrow \R$ such that the restrictions $u_{|(U \cup V)}=1$ and $u_{|(U_\varepsilon \cup V_\varepsilon)^c}=0$ hold and let finally $\mathcal{F}=\{g_1,g_2\}$ be the family given by Theorem \ref{SRM} for the manifold $M_2$. Define a map
 \begin{equation}\label{mapaefegorro}
 \hat{f}:\left( U_\varepsilon \cup V_\varepsilon \right) \times M_2 \rightarrow M_1 \times M_2  /  \hat{f}(x,y)=\left\{\begin{array}{c}
                                                            (F(x) , g_1(y)) \mbox{  if $x \in U_\varepsilon $} \\
                                                            (F(x) , g_2(y)) \mbox{  if $x \in V_\varepsilon $}
                                                          \end{array}\right.
 \end{equation}
and extend $\hat{f}$ to  \begin{equation}\label{mainmap}
f: M_1 \times M_2 \rightarrow M_1 \times M_2  /  f(x,y)=\left\{ \begin{array}{c}
   u(x).\hat{f}(x,y)+(1-u(x)).(F(x),y) \mbox{ if $x \in U_\varepsilon \cup V_\varepsilon$} \\
      (F(x),y) \mbox{\ \ \ \ \ \ \ \ \ \ \ \ \ \ \ \ \ \ \ \ \ \ \ \ \ \ \ \ \ \ \ \ \ \ \  if $x \notin U_\varepsilon \cup V_\varepsilon$.}
      \end{array} \right.
      \end{equation}
\subsection{Dynamics of $f$} In the next section we prove that $f$ is a $C^1$ robustly transitive map. The section begins highlighting a series of features the map $f$ possesses that are straightforward to check:

\begin{enumerate}\label{dynf}
  \item If $f(x,y)=(f_1(x,y),f_2(x,y))$, then $f_1(x,y)=F(x)$.
  \item The restriction $f_{|\left(U \cup V \right) \times M_2}=\hat{f}$.
  \item The restriction ${f_2}_{|U \times M_2}=g_1$.
  \item The restriction ${f_2}_{|V \times M_2}=g_2$.
\item The point $(p,r_1)$ is a source fixed point. Then, $f$ is not hyperbolic.
  \item The point $(p,a_1)$ is a saddle fixed point for $f$.
\item The local unstable set at $(p,a_1)$ is $W^u_{loc}(p,a_1)=U\times\{a_1\}$.
   \item The local stable set at $(p,a_1)$ is $W^s_{loc}(p,a_1)=\{p\}\times B$, where $B \subset M_2$ is an open ball containing $a_1$ and contained in the local stable manifold of $g_1$.
  \end{enumerate}

We prove ahead that both the local stable and local unstable sets of $(p,a_1)$ are dense in $M_1 \times M_2$. This will imply that the open set $U \times B$ has backward and forward iterates both dense yielding transitivity for $f$. Afterwards, since \textit{expanding} is a robust property and Theorem \ref{SRM} is robust we prove that $f$ is $C^1$ robustly transitive.

\begin{lema}\label{umd}
  The unstable set $W^u(p,a_1)$ is dense in $M_1 \times M_2$.
\end{lema}
\pf
Let $W=W_1 \times W_2$ be an open set in $M_1 \times M_2$. We show that there exists $j \in \N$ such that $f^j(W^u_{loc}(p,a_1)) \cap W \neq \emptyset$. Let $f$ be $f(x,y)=(F(x),f_2(x,y))$.\\ Since $F(U)=M_1$ then $f (W^u_{loc}(p,a_1))=f(U \times \{a_1\}) \supset M_1 \times \{a_1\}$.\\ Forward iteration by $f$ provides $f^{2}(W^u_{loc}(p,a_1)) \supset M_1 \times \{a_1,g_2(a_1)\}$ and also $f^{3}(W^u_{loc}(p,a_1)) \supset M_1 \times \{a_1,g_1g_2(a_1),g_2(a_1),g_2^2(a_1)\}$. It is straightforward noticing that $f^{k}(W^u_{loc}(p,a_1)) \supset M_1 \times \langle \mathcal{F}\rangle^k(a_1)$ for the $k$-th iterate of $a_1$ under the action of the semigroup $\langle \mathcal{F} \rangle ^+$. Now, since IFS$(\mathcal{F})$ is minimal, there exists a branch of the orbit of $ \langle \mathcal{F} \rangle^+ \left( a_1 \right)$ that intersects $W_2$ in a point $w_2$ after, say, $j$ iterates. Hence there exists $j \in \N$ such that $f^{j}(W^u_{loc}(p,a_1)) \supset M_1 \times \{w_2\}$ which implies $f^{j}(W^u_{loc}(p,a_1)) \cap W \neq \emptyset$.

\begin{lema}\label{smd}
  The stable set $W^s(p,a_1)$ is dense in $M_1 \times M_2$.
\end{lema}
\pf
Let $W=W_1 \times W_2$ be an open set in $M_1 \times M_2$ and $W^s_{loc}(p,a_1)=\{p\}\times B$ for some neighborhood $B$ of $a_1$ in $M_2$ contained in the stable manifold $W^s(a_1)$ of $g_1$ in $M_2$. We show that a backward iterate of $W^s_{loc}(p,a_1)$ intersects $W$.\\ Observe first that since $F$ is expanding, there exists a number $m \in \N$ satisfying $f^m(W)= F^m(W_1) \times f_2^m(W)=M_1 \times  f_2^m(W)$. Pick any point $x$ in $f_2^m(W)$. Since IFS$(\mathcal{F})$ is minimal, then $\langle \mathcal{F}\rangle^+(x)$ is dense in $M_2$, so there exists $n \in \N$ such that $\hat{f}_2^n(x) \cap B \neq \emptyset$ so $f^{m+n}(W)=M_1 \times f_2^n(f_2^{m}(W))$ satisfies $f^{m+n}(W)\cap \{p\}\times B \neq \emptyset$. Observe that since preimage of nonempty set might be empty, we can not deduce from the latter intersection that there is a point in $W^s_{loc}(p,a_1)$ with a backwards iterate in $W$. We construct such a point to finish the proof of the Lemma.\\ Consider the set $\bigcap_{t=0}^{t=n}f^{-t}(f^{m+n}(W))$, from item (2) at Remark \ref{preimagenes} follows that its projection onto $M_1 \cap (U \cup V)$ contains at least $2^{n+1}$ open and connected sets, obtained from a nested sequence of preimages of $M_1$ under $F$. One of these open and connected open sets, say $O$, realizes the itinerary of the branch where $W$ intersects $W^s_{loc}(p,a_1)$ going from $O \times f_2^{m}(W)$ up to $F^{n}(O) \times f_2^{m+n}(W)=M_1\times f_2^{m+n}(W)$. Notice that iteration under $f$ of $W$ still makes sense because being $F^m(W_1)=M_1$ it suffices to only consider the subset $W_0 \subset W_1$ satisfying $F^m(W_0)=O$ and reduce the starting set $W$ to its subset $W_0 \times W_2$ where everything until this point holds. Finally, since $F^m(W_0)=O$ and $F^n(O)=M_1$, there exist $o \in O$ and $x_1 \in W_0$ such that $F^n(o)=p$ and $F^m(x_1)=o$.\\ Summing up, there exists a point $(x_1,x_2) \in W_1\times W_2$ such that its iterates verify $f^{m+n}(x_1,x_2)=f^n(F^{m}(x_1),f_2^{m}(x_1, x_2))= (F^n(o),f_2^{m+n}(x_1, x_2)) \in \{p\}\times B $ which implies $f^{-m-n}(W^s_{loc}(p,a_1)) \cap W \neq \emptyset$.

\begin{rk}\label{nodependedeserabierto}
\begin{enumerate}\
  \item Lemma \ref{smd} holds for every set $W_1 \times W_2$ given $W_2$ be nonempty.
  \item The construction at the end of Lemma \ref{smd} gives a way to explicitly construct points satisfying the claim at the thesis of Lemma \ref{umd}.
  \item The set $U \times B$ is an open neighborhood of $(p,a_1)$ with dense forward and backward iterates.
\end{enumerate}
  \end{rk}

  \begin{thm}\label{fisRT}
  The map $f$ defined by Equation (\ref{mainmap}) is $C^1$ robustly transitive.
\end{thm}
\pf According to Lemmas \ref{umd} and \ref{smd} there exists a point in $M_1 \times M_2$ with an open neighborhood dense under forward and backward iteration. Hence it holds that for all open sets $A$ and $B$ in $M_1 \times M_2$ there exists a natural number $k$ such that $f^k(A) \cap B \neq \emptyset$ which yields transitivity for $f$.\\ To prove robustness, consider any $C^1$ neighborhood $\mathcal{V}_{f}$ of $f$. Observe that for the map $\hat{f}$ defined by Equation (\ref{mapaefegorro}), by Theorem \ref{SRM} there exists a $C^1$ neighborhood $\mathcal{V}_{\hat{f}}$ of $\hat{f}$ such that every $\hat{g}$ in $\mathcal{V}_{\hat{f}}$ induces a map $g$ defined as by Equation (\ref{mapaefe}) with $g$ satisfying Lemmas \ref{umd} and \ref{smd}, so there is $C^1$ robustness of the construction that defines $f$ from $\hat{f}$. Thereafter, take $\mathcal{V}_{f}$ and reduce it to a neighborhood $\mathcal{W}_{f}$ of $f$ such that every map in $\mathcal{W}_{f}$ is, restricted to $(U \cup V) \times M_2$, contained in $\mathcal{V}_{\hat{f}}$ (this is possible after item 2 at the beginning of Subsection \ref{dynf}).
Finally, since $f=(F,f_2)$ and there exists a neighborhood $\mathcal{U}_F$ of $F$ in the $C^1$ topology containing only expanding maps \cite[Corollary to Theorem ($\alpha$)]{shub}, reduce $\mathcal{W}_{f}$ to a neighborhood $\mathcal{U}$ until all maps in $\mathcal{U}$ are expanding when projected to $M_1$. In fact, all we need to hold is that every map $h$ in $\mathcal{U}$ satisfies $h(U \times M_2)=h(V \times M_2)=M_1 \times M_2$ which is almost straightforward by taking a power of $F$, if needed, at the beginning of the construction. In turn, we have a $C^1$ neighborhood $\mathcal{U}$ of $f$ such that all of its maps satisfy Lemmas \ref{umd} and \ref{smd} which gives robust transitivity for $f$.

\section{A singular endomorphism $A$ of $M_1 \times M_2$.}
 We procceed now to the second step of the construction: perform a perturbation on the map $f$ defined by Equation (\ref{mainmap}) in order to endow it with robust singularities without destroying its dynamical characteristics. Recall from the introduction the definitions of critical point and critical set. Recall also that a map $h$ is a \textit{robustly singular} endomorphism if there exists a neighborhood $\mathcal{U}_h$ of $h$ in the $C^1$ topology such that all $ g \in \mathcal{U}_h$ satisfy $S_g \neq \emptyset$.

\subsection{Construction of $A$}\label{constructionofF} Let $s$ be a point not in $\left( U_\varepsilon \cup V_\varepsilon \right) \times M_2$. Fix a small ball around $s$ and perturb $f$ to a map $A$ by introducing critical points artificially inside the ball in a robust way. The perturbation is such that it does not affect $f$ on the first factor neither on the blending region so $A$ inherits the robust transitivity of $f$.\\
Before going into the details of the construction we want to point out to the reader that the construction carried on in \cite[Section 4]{mo2} is, with some adjustments, the same that we perform to produce the robustly singular map perturbed from $f$. Our approach will be to take the point $s$ in our manifold, take a chart so that we can work on a ball in $M_1 \times M_2$ through coordinates in $\R^{m_1+m_2}$ (it is here when we use the fact that $M_2$ is a real manifold), and perform the surgery to introduce the robust singularities. The reader is invited to go over \cite[Section 4]{mo2} to get in touch with a simpler version of the construction ahead.\\

Let $\mu:M_1 \times M_2 \to \R^{m_1+m_2}$ be a chart such that $\mu:=(\mu_1,\mu_2)$, satisfying that $\mu_1 \left( U_\varepsilon \cup V_\varepsilon \right) \subset B_{\left(0,\frac{1}{10}\right)}$ and $\mu_1(p)=0$. Let $s:=(s_1, s_2) \in M_1 \times M_2$ be such that $\mu(s)=(\frac{1}{4},0,...,0,\frac{1}{4})$. Since $||\mu_1 (s_1)||=\frac 14$ then $\mu_1(s_1) \notin B_{\left(0,\frac{1}{10}\right)}$ so $s$ does not belong to the blending region $\left( U_\varepsilon \cup V_\varepsilon \right) \times M_2$ either.\\ For the rest of the section we adopt the following extended notation for our chart: $\mu(x,y)=(\mu_1(x),\mu_2(y))=(\mu_{11}(x),...,\mu_{1m_{1}}(x),\mu_{21}(y),...,\mu_{2m_{2}}(y))$.

\begin{figure}[ht]
\begin{center}
\subfigure[]{\includegraphics[scale=1.5]{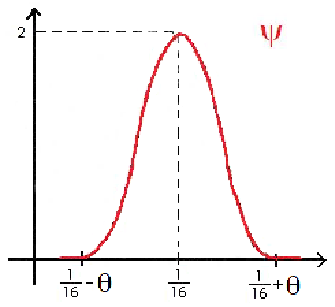}}
\subfigure[]{\includegraphics[scale=0.6]{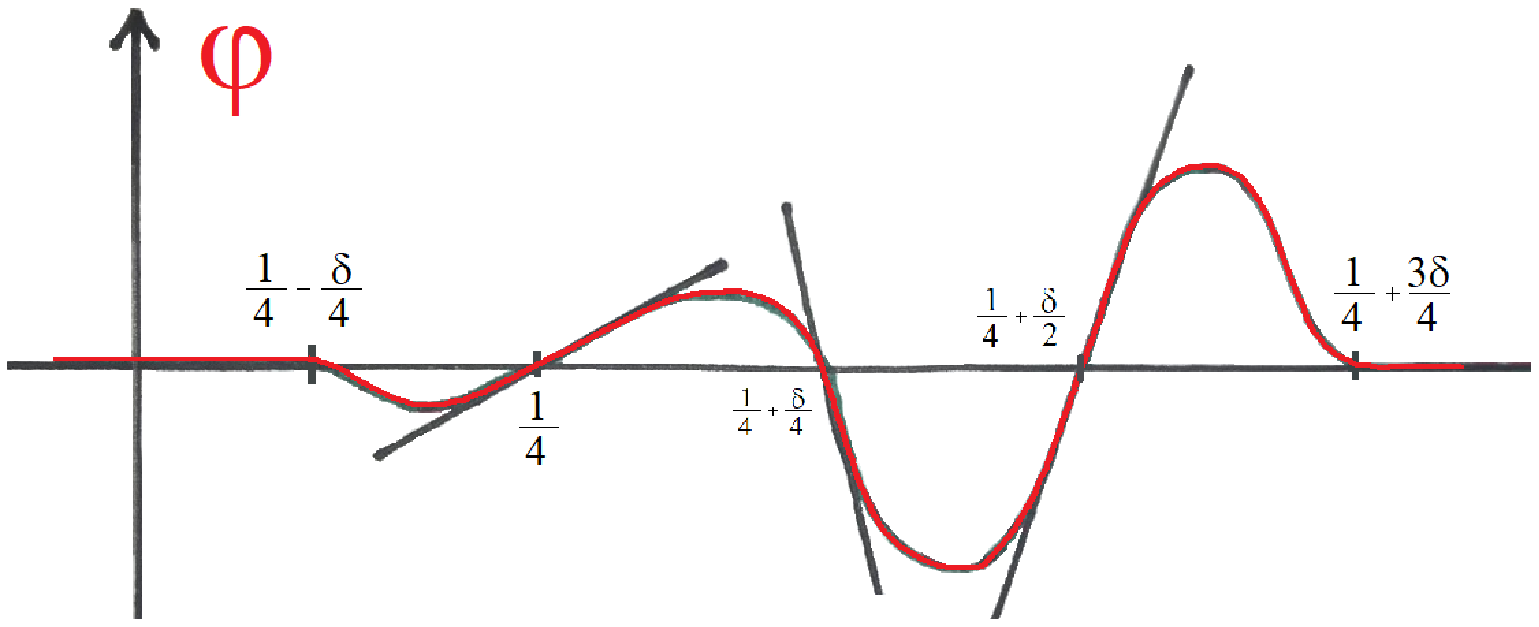}}
\caption{Graphs of $\psi$ and $\varphi$ }\label{figura11}
\end{center}
\end{figure}

We set now a series of (technical) parameters and considerations that lead to the definition of $A$. Start fixing $r>0$ satisfying $B_{(s,r)} \cap [\left( U_\varepsilon \cup V_\varepsilon \right) \times M_2] = \emptyset$, this is possible since $s \notin  (U_\varepsilon \cup V_\varepsilon) \times M_2 $. Observe that $B_{(\mu_1(s_1),r)} \cap B_{(0,\frac {1}{10})} = \emptyset$. Fix a second parameter $\theta$ such that  $0< 2\theta  < r$  and finish fixing a third parameter $\delta$ such that $0<\delta < 2\theta $. Observe that for the points $q_1, q_2 \in M_2$ satisfying that $\mu_2(q_1)=(0,0,..,0, \frac 14 + \frac \delta 4)$ and $\mu_2(q_2)=(0,0,..,0, \frac 14 + \frac \delta 2)$ holds that both $(s_1,q_1)$ and $(s_1,q_2)$ belong to $B_{(s,r)}$. \\
Define now a smooth function $\psi :\R\to\R$ with a unique critical point at $\frac{1}{16}$, with $\psi (\frac{1}{16})=2$ and $\psi (t)=0$ for all $ t$ in the complement of $(\frac{1}{16}-\theta ,\frac{1}{16}+\theta)$ as shown in Figure \ref{figura11} (a), and another smooth function $\varphi:\R\to\R$ such that:
\begin{itemize}
\item $\varphi $ is as in Figure \ref{figura11} (b),
\item $\varphi (\frac{1}{4})=\varphi (\frac 14 +\frac \delta 4)=\varphi (\frac 14 +\frac \delta 2)=0$,
\item $\varphi' (\frac{1}{4})=\frac{1}{2}$, $\varphi' (\frac{1}{4} +\frac{\delta}{4})=-3$, $\varphi' (\frac{1}{4} +\frac{\delta}{2}  )=3$,
\item $\varphi (t)=0$ for all $ t \notin [\frac{1}{4}-\frac{\delta}{4},\frac{1}{4}+\frac{3\delta}{4} ]$.
\end{itemize}

We have all that is needed to define a perturbation of $f$ on its last entry, depending on the parameters $r$,$\theta$ and $\delta$. A new map $A(r,\theta,\delta): M_1 \times M_2 \to M_1 \times M_2$ denoted just by $A$ will be obtained, defined at the point $(x,y)=(x_1,...,x_{m_1},y_1,...,y_{m_2})$ by

\begin{equation}\label{mapaefe}
            A(x,y)=\left\{\begin{array}{c}
                   f(x,y) \mbox{\ \ \ \ \ \ \ \ \ \ \ \ \ \ \ \  \ \ \ \ \ \ \ \ \ \ \ \ \ \ \ \ \ \ \ \ \ \ \ \ \ \ \ \ \ \ \ \ \ \ \ \ \ \ \ \ \ \ \ \ \ \ \ \ \ \ \ \ \ \ \ \ \ \ \ \  if $(x,y) \notin B_{(s,r)}$} \\
                   \left(F(x),y_1,...,y_{m_2-1},\mu_{2m_2}^{-1}\left( \mu_{2m_2}(y) - {\varphi(\mu_{2m_2}(y))}.\psi ( ||\mu_1(x)||^2 )\right)\right) \mbox{ if $(x,y) \in B_{(s,r)}$}
                 \end{array} \right.
         \end{equation}

\begin{rk}\label{rkA}
\begin{enumerate}\
  \item For all $(x,y) \notin (U_\varepsilon \cup V_\varepsilon)\times M_2$ it holds that $f(x,y)=(F(x),y)$.
   \item For all $(x,y) \notin B_{\left(s,\frac{3\delta}{4}\right)}$ it holds that $A(x,y)=f(x,y)$.
   \item If $(x,y)$ in $M_1 \times M_2$ and $A(x,y)=(A_1(x,y),A_2(x,y))$, then $A_1(x,y)=F(x)$. Therefore it exists a $C^1$ neighborhood $\mathcal{U}_A$ of $A$ such that for all $h \in \mathcal{U}_A$ it holds that $h(U \times M_2)=h(V \times M_2)=M_1 \times M_2$ and for all $W$ open in $M_1 \times M_2$ it holds that the projection of $h(W)$ to $M_1$ has nonempty interior in $M_1$.
      \item In the rest of the section we denote $\varphi(\mu_{2m_2}(y))$ as $\varphi$, and $\psi\left( ||\mu_1(x)||^2 \right)$ as $\psi$, omitting the evaluations appearing on the definition to make the reading easier.

\end{enumerate}
\end{rk}

\begin{lema}\label{Fps}
The map $A$ defined by Equation (\ref{mapaefe}) is robustly singular.
\end{lema}
\pf
Start computing the differential $D_{(x,y)} A$ at all points $(x,y) \in B_{(s,r)} $. By Equation (\ref{mapaefe}) we have that, after taking determinants, it holds
$$det(D_{(x,y)}A) =det(D_xF).\left(\frac{\partial \mu_{2m_2}^{-1}}{\partial y_{m_2}}(\mu_{2m_2}(y)-\varphi.\psi)\right).\left(\frac{\partial \mu_{2m_2}}{\partial y_{m_2}}(y)\right) .\left( 1-\varphi'.\psi \right).$$

Perform calculations at points where $\varphi =0$ and $\varphi' \neq 0$ applying the chain rule. We have that $det(D_{(s_1,q_1)}A)=7.det(D_{s_1}F)$ and $det(D_{(s_1,q_2)}A)=-5.det(D_{s_1}F)$ as well as $det(D_{(s_1,s_2)}A)=0$. Since $D_{s_1}(F)>1$ there exists a $C^1$ neighborhood $\mathcal{U}_A$ of $A$ of radio $1$ such that for every $ g \in \mathcal{U}_A$, $S_g$ is nonempty so $A$ is robustly singular.

\subsection{Dynamics of $A$.}\label{dynA} We finish our construction with the proof of robust transitivity for $A$ and the proof of Theorem \ref{main}.

\begin{lema}\label{Frt}
  The map $A$ defined by Equation (\ref{mapaefe}) is $C^1$ robustly transitive.
\end{lema}
\pf By Remark \ref{rkA} and Remark \ref{nodependedeserabierto} follow that Theorem \ref{fisRT} holds for $A$. Furthermore, choose any $C^1$ neighborhood of $A$ and reduce it exactly as it was done in Theorem \ref{fisRT} to get a neighborhood satisfying the claim at the thesis of Lemma \ref{Frt}.

\begin{proof}[ \textbf{Proof of Theorem \ref{main}}]
  Let $A$ be the endomorphism defined by Equation (\ref{mapaefe}). Define $\mathcal{U}_1$ a $C^1$ open neighborhood of $A$ where Lemma \ref{Fps} holds and $\mathcal{U}_2$ a $C^1$ open neighborhood of $A$ where Lemma \ref{Frt} holds. Then, all maps in $\mathcal{U}_A=\mathcal{U}_1 \cap \mathcal{U}_2$ are transitive and have nonempty critical set, hence $A$ is a robustly singular and $C^1$ robustly transitive endomorphism supported on $M_1 \times M_2$.
\end{proof}

\section*{Acknowledgements} The author would like to give thanks to Prof. Jorge Iglesias for fruitful conversations regarding the problem addressed by this article. The author would also like to thank Prof. Roberto Markarian and Prof. Lorenzo J. D\'iaz for their generous and encouraging attitude towards the work of the author.

\end{document}